\documentclass[12pt,twoside]{amsart}

\usepackage{amssymb,amsmath}
\theoremstyle{plain}
\newtheorem{theorem}{Theorem}[section]
\newtheorem{lemma}[theorem]{Lemma}

\theoremstyle{definition}

\newtheorem{defn}[theorem]{Definition}

\numberwithin{equation}{section}

  \input xypic
  \xyoption{all}

 \usepackage{tikz}
\usetikzlibrary{arrows}

 \usepackage{times}

  \usepackage{amsmath, amsthm, amsfonts}

\def \df{\square}
\def \mc{\mathcal}
\def \inv{^{-1}}

\def \cplane{\mathbb{C}}
\def \mff{\mathsf}

\def \al{\frac{1}{3}}
\def \six{\frac{1}{6}}
\def \half{\frac{1}{2}}

\def \b{\bar }

\def \p{\partial}

 \def\CHART{\mathsf{U}}

\def \p{\partial}

\def \v {\vskip 0.1in}
\def \n {\noindent}

\def \df{\square}
\def \mc{\mathcal}
\def \inv{^{-1}}

\def \v{\vskip 0.1in}
\def \n{\noindent}

\def \cplane{\mathbb{C}}

\def \al{\frac{1}{3}}
\def \six{\frac{1}{6}}
\def \half{\frac{1}{2}}

\def \b{\bar }

\def \p{\partial}

 \def\CHART{\mathsf{U}}

\begin{document}

\title{Differential inequalities on homogeneous toric bundles}

\author[Li]{An-Min Li}
\address{Department of Mathematics,
Sichuan University,
 Chengdu, 610064, China}
\email{anminliscu@126.com}

\author[Sheng]{Li Sheng  }
\address[Sheng]{Department of Mathematics Sichuan University Chengdu, 610064, China}
\email[Corresponding author]{lshengscu@gmail.com}

\author[Zhao]{Guosong Zhao}
\address[Zhao]{Department of Mathematics Sichuan University
Chengdu, 610064, China}
\email{gszhao@scu.edu.cn}
\thanks{ Li acknowledges the support of NSFC Grant 11521061. \\
${}\quad\ $Sheng acknowledges the support of NSFC Grant 11471225.\\
${}\quad\ $Zhao acknowledges the support of NSFC Grants 11571242.
}

 {\abstract
We study a generalized Abreu Equation in $n$-dimensional polytopes and prove some differential inequalities for
homogeneous toric bundles.
\endabstract}

\maketitle
\n



\section{Introduction}\label{Sec-Intro}
The existence of extremal and contant scalar curvature is a central problem in K\"ahler geometry. In a series of papers \cite{D1}, \cite{D2}, \cite{D3}, and \cite{D4}, Donaldson studied
this problem on toric manifolds and proved the existence
of metrics of constant scaler curvatures on toric surfaces
under an appropriate stability condition. Later on in \cite{CLS1} and \cite{CLS2}, Chen, Li and Sheng proved the existence of metrics of prescribed scaler curvatures on toric surfaces under the uniform stability condition.
\v
It is important to generalize the results of Chen, Li and Sheng to more general K\"ahler manifold. This is one of a sequence of papers, aiming at generalizing the results of Chen, Li and Sheng to homogeneous toric bundles.
 A homogeneous toric bundle $G\times_{K} M$ have a compact toric manifold $M$ as fiber and a generalized flag manifold $G/K$ as basis. Let $U\subset M$ be any complex coordinate domain and let $G(U)\subset G\times_{K} M$ be the image of $U$ under the left $G$-action on $G\times_{K} M$.
In this paper we establish some differential inequalities on $G(U)$, which are generalizations of the differential inequalities in \cite{CLS1}. The differential inequalities of this paper play a important role in our works on scaler curvatures for homogeneous toric bundles.

\section{Homogeneous toric bundles}\label{sec-Homogeneous}

We recall some facts about homogeneous toric bundles and refer to \cite{PS}, \cite{Ar} and \cite{D5} for details.
Let $G$ be a compact semisimple Lie group, $K$ be the centralizer of a torus $S$ in $G$. Let $T$ be a maximal torus in $G$ containing $S$. Then $T\subset C(S)=K$, $G/K$ is a generalized flag manifold. Denote $\mathfrak{o}=eK$.  Let  $\mathfrak{g}$ (resp.  $\mathfrak{k}$, $\mathfrak{h}$ ) be the Lie algebra of $G$ (resp. $K$, $T$). Let $B$ denote the Killing form of $\mathfrak{g}$. Recall that $-B$ is a positive definite inner product on $\mathfrak{g}$. There is a orthogonal decomposition with respect to $-B$:
$$\mathfrak{g}=\mathfrak{k}\oplus \mathfrak{m},\;\;\;Ad(k)\mathfrak{m}\subset \mathfrak{m}\;\;for\;all\;k\in K.$$
The tangent space of $G/K$ at $\mathfrak{o}$ is identified with $\mathfrak{m}$.
\v
Let $Z(K)$ be the center of $K$, it is an n-dimensional torus, denoted by $T^n$. Let $(M,\omega)$ be a compact toric K\"{a}hler manifold of complex dimension $n$, where $T^n$ acts effectively on $M$. Let $\varrho: K\rightarrow T^n$ be a surjective homomorphism. The homogeneous toric bundle $G\times_{K}M$ is defined to be the space $G\times M$ modulo the relation $(gh,x)=(g,\varrho(h)x)$,
where $g\in G$, $h\in K$, and $x\in M$. Later we will omit $\varrho$ to simplify notations. The space $G\times_{K}M$ is a fiber bundle with fiber $M$ and base space $G/K$, a generalized flag manifold. There is a natural $G-$action on $G\times_{K}M$ given by $g\cdot[h,x]=[gh,x]$, $g\in G$, $x\in M$, and a natural $T^n-$action on $G\times_{K}M$ given by $k\cdot[g,x]=[g, k^{-1}x]$, $k\in T^n$. Both $G/K$ and $G\times_{K}M$ are complex manifolds. In fact, they can be expressed as $G^{\mathbb{C}}/P$ and $G^{\mathbb{C}}\times_P M$ respectively, where $G^{\mathbb{C}}$ is the complexification of $G$, $P$ is a parabolic subgroup of $G^{\mathbb{C}}$.
\v

Denote by $\mathfrak{g}^{\mathbb {C}},\mathfrak{k}^{\mathbb {C}}, \mathfrak{h}^{\mathbb {C}}$ the complexification of   $\mathfrak{g}$, $\mathfrak{k}$, $\mathfrak{h}$ respectively.
Let $R$ be the root system of $\mathfrak{g}^{\mathbb {C}}$ with respect to $\mathfrak{h}^{\mathbb {C}}$, we have the root space decomposition
$$\mathfrak{g}^{\mathbb {C}}=\mathfrak{h}^{\mathbb {C}}\oplus \sum_{\alpha\in R}\mathbb {C}E_{\alpha}=\mathfrak{k}^{\mathbb {C}}\oplus \mathfrak{m}^{\mathbb {C}},$$
$$\mathfrak{k}^{\mathbb {C}}=\mathfrak{h}^{\mathbb {C}} \oplus \sum_{\alpha\in R_{K}}\mathbb {C}E_{\alpha},\;\;\; \mathfrak{m}^{\mathbb {C}}=\sum_{\alpha\in R_{M}}\mathbb {C}E_{\alpha},$$
where $R_{K}$ is a subset of $R$ and $R_{M}=R\setminus R_{K}$. For any $\varphi\in \mathfrak{h}^{*}$ we define $h_{\varphi}\in \mathfrak{h}$ by
$$B(h, h_{\varphi})=\varphi(h)\;\;\;for\;\;all\;\;h\in \mathfrak{h},$$
and set
$$H_{\varphi}=\sqrt{-1}h_{\varphi}.$$
Define
$$\mathfrak{t}=\mathfrak{h}\bigcap \mathfrak{z}(\mathfrak{k}^{\mathbb{C}}),$$
Let
$$\kappa:\mathfrak{h}^*\rightarrow \mathfrak{t}^* \;\;\;\alpha \mapsto \alpha|_{\mathfrak{t}}$$
be the restriction map. Set $R_T=\kappa(R)=\kappa(R_{M})$. The elements of $R_T$ are called $T$-roots.

Choose a fundamental system $\Pi=\{\alpha_1,\cdots,\alpha_{n},\alpha_{n+1},\cdots,\alpha_{n+\hbar}\}$ of $R$. We fix a lexicographic ordering and let $R^+$ (resp. $R_K^+$)  be the set of positive roots of $R$ ( resp. $R_K$ ) with respect to $\Pi$. Let $R^+_M=R^+\setminus R_K^+$.
Denote $\Pi_{K}=\Pi\cap R_{K}$. We assume that $\Pi_{K}=\{\alpha_{n+1},\cdots,\alpha_{n+\hbar}\}.$
For $i\leq n,$ let
$$\tilde{h}_{\alpha_i}= h_{\alpha_i}+\sum_{j=n+1}^{n+\hbar} d_{i}^{j}h_{\alpha_j},$$ where $d_{i}^{j}$ are constants such that $B( \tilde{h}_{\alpha_i},h_{\alpha_j})=0$ for all $j=n+1,\cdots, n+\hbar.$ Denote
$$
H_{\alpha_i}= {\sqrt{-1}}h_{\alpha_{i}},\;\;i=n+1,\cdots, n+\hbar,\;\;\; \tilde{H}_{\alpha_j}= {\sqrt{-1}}\tilde{h}_{\alpha_{j}},\;\;j=1,\cdots, n.
$$
For any $\alpha\in R_{K}$, we have $\alpha(\tilde{H}_{\alpha_j})=0.$ Then  $\mathfrak{t}=span\{\tilde{ H}_{\alpha_1},\cdots, \tilde{H}_{\alpha_n}\}.$
Let $\mathfrak{h}'=span\{ H_{\alpha_{n+1}},\cdots, H_{\alpha_{n+\hbar}}\}.$ We have the orthogonal decomposition with respect to $-B$:
$\mathfrak{h}=\mathfrak{t}\oplus \mathfrak{h}'.$
\v

\v
We choose a Weyl basis $e_{\alpha}\in \mathfrak{g}_{\alpha}^{\mathbb{C}}$ of $\mathfrak{g}^{\mathbb{C}}$
such that, for $e_{\alpha}\in \mathfrak{g}_{\alpha}$ and $e_{-\alpha}\in \mathfrak{g}_{-\alpha}$,
\begin{equation}\label{eqn_2.1}
B(e_{\alpha},e_{-\alpha})=1,\;\;\;[e_{\alpha},e_{-\alpha}]=h_{\alpha},\end{equation}

Obviously, $[h,e_{\alpha}]=\alpha(h)e_{\alpha},$ for any $h\in \mathfrak{t}$. Set  $$  V_{\alpha}=e_{\alpha}-e_{-\alpha}, \;\;\;\;\; {W}_{\alpha}=\sqrt{-1}(e_{\alpha}+e_{-\alpha}),\;\;\;\;\;
  H_{\alpha}= \sqrt{-1}h_{\alpha}.$$
It is easy to see that $  H_{\alpha},  V_{\alpha},  {W}_{\alpha} \in (\mathfrak{g}_{\alpha}\oplus \mathfrak{g}_{-\alpha}\oplus [\mathfrak{g}_{\alpha},\mathfrak{g}_{-\alpha}])\cap \mathfrak{g}.$
For any $\alpha\in  {R}^{+}$ we have  $\alpha= \sum_{j=1}^{n+\hbar} M_{\alpha}^{j}\alpha_{j}$ with $M_{\alpha}^{j}\geq 0.$ It is easy to see that
\begin{equation}\label{eqn_4.6}
H_{\alpha}= \sum_{j=1}^{n+\hbar} M_{\alpha}^{j}H_{\alpha_j}= \sum_{j=1}^{n} M_{\alpha}^{j}\tilde{H}_{\alpha_j}+ \sum_{j=n+1}^{n+\hbar}M_{\alpha}^{'j}H_{\alpha_j},
\end{equation}
where $M_{\alpha}^{'j}=M_{\alpha}^{j}-\sum_{i=1}^n M_{\alpha}^{i}d_i^j$.
Obviously, for any $\alpha\in R_{M^{+}}$, we have
\begin{equation}
\sum_{j=1}^{n} M_{\alpha}^{j}>0.
\end{equation}
\v\v

\section{$(G,T^n)$-invariant K\"ahler metrics}\label{sec-Kahler metrics}
\v
From now on, our convention for the ranges of indices is the following:
$$1\leq A,B,...\leq n+l $$
$$1\leq i,j,k,...\leq n,$$ $$ n+1\leq \alpha \leq n+l$$
where $l $ is the dimension of $\mathfrak{m}.$
For any $1\leq j\leq n,\alpha\in R_M^{+}$, let $\tilde{H}^*_{\alpha_j} $, $  V^*_{\alpha},  W^*_{\alpha}   $ be the fundamental vector fields  corresponding to  $\tilde{H}_{\alpha_j} $, $V_{\alpha}, W_{\alpha}$.
Then the left-invariant vector fields $\{\frac{\p}{\p x^{j}}, \tilde{H}^*_{\alpha_j},V^*_{\alpha}, W^*_{\alpha}\}_{1\leq j\leq n,\alpha\in R_M^{+}}$ is a local base of $G\times_{K}M$. Obviously,
\begin{equation}
[\tfrac{\p}{\p x^{j}}, \tilde{H}^*_{j}]=0,\;\;\; [\tfrac{\p}{\p x^{j}}, V^*_{\alpha}]=[\tfrac{\p}{\p x^{j}}, W^*_{\alpha}]=0,\;\;\;1\leq j\leq n,\alpha\in R_{M^{+}}.
\end{equation}
Let $\{dx^{j},\nu^{j}, dV^{\alpha},dW^{\alpha}\}_{1\leq j\leq n,\alpha\in  R_M^{+}}$ be the dual left-invariant 1-form of the base.
Set
$$S_{j}=\frac{1}{2}(\frac{\p}{\p x^{j}}-\sqrt{-1} H_{j}^*),\;\;\;\;\; 1\leq j\leq n$$
$$S_{\alpha}=\frac{1}{2}(  V^*_{\alpha}-\sqrt{-1}  W^*_{\alpha}),\;\;\;\;\alpha\in R_{M^{+}}.$$ We define the almost complex structure $J$ by,
$$
JS_{i}=\sqrt{-1}S_{i},\;\;\; JS_{\alpha}=\sqrt{-1}S_{\alpha},\;\;\;\;\forall\; 1\leq i\leq n,\;\;\alpha\in R_{M^{+}}.
$$
Then $J$ is a $(G,T^n)$-invariant complex structure. It
is easy to see that $S_{j}, S_{\alpha}$ are $(1,0)$-vector fields. Denote $\omega^{j},\omega^{\alpha}$ the dual $(1,0)$-form of $S_{j}, S_{\alpha}$. It is easy to see that  $S_{\alpha}$, $\bar S_{\alpha}$ is induced vector field of $e_{\alpha}$ and $-e_{-\alpha}.$
  By \eqref{eqn_2.1} we have
 \begin{align}
 \label{eqn_c_4.7}& [S_{j},\bar S_{k}]= [S_{j}, S_{k}]=0,\;\;\;[S_{\alpha},\bar S_{\alpha}]=-\sqrt{-1}H_{\alpha}^{*}=\sum_{j} M_{\alpha}^{j}(S_{j}-\bar S_{j}),\\
 &[S_{j},S_{\alpha}]=-\tfrac{1}{2}\alpha(\tilde h_{j})S_{\alpha},\;\;\;[S_{j},\bar S_{\alpha}]=\tfrac{1}{2}\alpha(\tilde h_{j})\bar S_{\alpha},  \;\;\; [S_{\alpha},S_{\beta}]=-N_{\alpha\beta}S_{\alpha+\beta},
\label{eqn_c_4.8}  \\
 &[S_{\alpha},\bar S_{\beta}]=N_{\alpha-\beta}S_{\alpha-\beta},\;\;\alpha>\beta\;\;\;   [S_{\alpha}, \bar S_{\beta}]=-N_{\alpha-\beta} \bar S_{\beta-\alpha}, \;\;\;\; \alpha<\beta. \label{eqn_c_4.9}
 \end{align}
As in \cite{N-1} one can check that
\begin{equation}\label{eqn_nu_1}
d\nu^{j}=\sum_{\alpha\in R_{M^{+}}}2 M_{\alpha}^{j}dV^{\alpha}\wedge dW^{\alpha}=\sum_{\alpha\in R_{M^{+}}}M_{\alpha}^{j}\sqrt{-1}\omega^{\alpha}\wedge \bar\omega^{\alpha}.
\end{equation}
In fact,   we have $$d\nu^{j}(S_{A},\bar S_{B})=S_{A}(\nu^{j}(\bar S_{B}))-\bar S_{B}(\nu^{j}(S_{A}))-\nu^{j}([S_{A},\bar S_{B}])=-\nu^{j}([S_{A},\bar S_{B}]).$$ Then \eqref{eqn_nu_1} follows from \eqref{eqn_c_4.7}-\eqref{eqn_c_4.9}.
\v
Denote by $\tau:M\rightarrow \bar{\Delta}\subset \mathfrak{t}^*$
the moment map of $M$, where $\Delta$ is a Delzant polytope.
The left invariant 1-form $\{\nu^{1},\cdots, \nu^{n}\}$ can be seen as a basis of $\mathfrak{t}^{*}.$
The moment map $\tau: M \rightarrow  \mathfrak{t}^{*}$
has components, relative to this basis of $\mathfrak{t}^{*}$, which we denote by $\tau_i$. Note that $\sum_{i=1}^{n} \tau_{i}\nu^{i}$ is independent of the choice of the bases.
\v
Now we fit $G$ into this picture. We fix a point $o\in \mathbb{R}^n$ and
identify $\mathfrak{t}^*$ with $T_{o}\mathbb{R}^{n}$. We choose $\{o,\nu^{i},\;i=1,...,n\}$ as a base of $\mathbb{R}^n$.
Let $\xi=(\xi_1,...,\xi_n)$ be the coordinate system with respect to the bases.
We choose $  \bar{\Delta}\subset \{(\xi_1,...,\xi_n)| \xi_1>0,\;\xi_2>0,\;...,\xi_n>0\}$ such that
\begin{equation}\label{equ_R_3.6}
 \sum_{\alpha\in R_{M^{+}}}\frac{\sum_{j=1}^{n} M_{\alpha}^{j}diam(\Delta)}{D_{\alpha}}<\frac{n}{4},
\end{equation}
where $$D_{\alpha}:=2\sum_{i=1}^{n} \tau_{i}\nu^{i}=2\sum_{j=1}^{n} M_{\alpha}^{j}\xi_{j}>0\;\;\forall \;\xi\in \bar{\Delta}.$$
Since the moment map is equivariant we can also regard $\tau$
as a map from $G\times_{K}M$ to $\mathfrak{t}^*$ and the components $\tau_{i}$ as functions on $G\times_KM$, that is, we extend $\tau: G\times_{K}M\rightarrow \bar{\Delta}$ by $\tau([g,x])=\tau(x)$. Following Donaldson (\cite{D4}) we consider the following form
\begin{equation*}
\Omega=d\left(\sum_{i=1}^{n} \tau_{i}\nu^{i}\right)
\end{equation*}
in $\tau^{-1}(\Delta).$
By \eqref{eqn_nu_1} we can write it as
\begin{equation}\label{eqn_4_sym}
\Omega=\sum_{i=1}^{n} d\tau_{i}\wedge\nu^{i} +\sum_{\alpha\in R_{M^{+}}} D_{\alpha}(dV^{\alpha}\wedge dW^{\alpha}).
\end{equation}
 It is easy to see that
$J$ is $\Omega$-compatible. i.e.,
$$
 \Omega(X,J{X})> 0, \;\;\;\forall X\neq 0,\;\;\Omega(JX,JY)=\Omega(X,Y).
 $$
So $(G\times_{K}M, \Omega)$ is a K\"ahler manifold with the K\"ahler form $\Omega$.
It is well known that there exists a convex function $f$ such that
$\frac{\p ^2 f}{\p x^{i}\p x^{j}}=\frac{\p \tau_{i}}{\p x^{j}}$. Then the K\"ahler form $\Omega$ can be written as
\begin{align}
\Omega&=\sum_{i,j=1}^{n} f_{ij}dx^{i}\wedge \nu^{j} +\sum_{\alpha\in R_{M^{+}}}D_{\alpha}(dV^{\alpha}\wedge dW^{\alpha}) \nonumber\\
&= \frac{1}{2}\left[\sum_{i,j=1}^{n}f_{ij}\omega^{i}\wedge \bar \omega^{j} +\sum_{\alpha\in R_{M^{+}}}D_{\alpha}(\omega^{\alpha}\wedge \bar \omega^{\alpha})\right].
\end{align}
We denote $\Omega$ by $\Omega_{f}$. Then the Riemannian metric is given by  \begin{align*}
\mathcal G_{{f}}&= \sum_{i,j=1}^{n} f_{ij}(dx^{i}\otimes dx^{j} +\nu^{i}\otimes \nu^{j}) +\sum_{\alpha\in R_{M^{+}}} D_{\alpha}(dV^{\alpha}\otimes dV^{\alpha}+dW^{\alpha}\otimes dW^{\alpha})\\
&=\sum_{i,j} \frac{f_{ij}}{2}(\omega^{i}\otimes \bar \omega^{j}+\bar \omega^{j}\otimes \omega^{i})+\sum_{\alpha\in R_{M^{+}}} \frac{D_{\alpha}}{2}\left(\omega^{\alpha}\otimes   \bar \omega^{\alpha}+  \bar \omega^{\alpha} \otimes \omega^{\alpha}\right).
\end{align*}
Let
$$\xi_i=\frac{\p f}{\p x^{i}},\;\;u(\xi_1,...,\xi_n)=\sum_{i=1}^n x^{i}\xi_i - f(x).$$
Then
 \begin{align*}
\mathcal G_{u}&= \sum_{i,j=1}^{n}( u_{ij}d\xi_{i}\otimes d\xi_{j} + u^{ij}\nu^{i}\otimes \nu^{j}) +\sum_{\alpha\in R_{M^{+}}} D_{\alpha}(dV^{\alpha}\otimes dV^{\alpha}+dW^{\alpha}\otimes dW^{\alpha}),
\end{align*}
where $u_{ij}=\frac{\p^2 u}{\p \xi_i \p \xi_j} $ and $(u^{ij})$ is the inverse matrix of $(u_{ij})$.
The K\"ahler form $\Omega$ and the K\"ahler metric $\mathcal G$ can be extended over $G\times_{K}M$.
\v
Suppose that $\Delta$ is defined by linear inequalities $h_k(\xi)-c_k>0$,  for $k=1, \cdots, d$,
where $c_k$ are constants and
$h_k$ are affine linear functions in $\mathbb  R^n$, $k=1, \cdots, d$, and each $h_k(\xi)-c_k=0$ defines a facet of $\Delta$.
Write $\delta_k(\xi)=h_k(\xi)-c_k$
and set
\begin{equation}\label{eqn2.1}
v(\xi)=\sum_k\delta_k(\xi)\log\delta_k(\xi).
\end{equation}
It defines a K\"ahler metric on $G\times_{K}M$, which we call the {Guillemin} metric.

\v\v
\section{Generalized Abreu Equations}\label{sec-Abreu Equations}

In this section we calculate the Ricci tensor and the scaler curvature of $G\times_KM$. Podesta and Spiro have calculated the Ricci tensor in \cite{PS}. Here we use a different method. First we note that
for any $(1,0)$-vector fields $X,Y$, we have $\nabla_{X}Y $ is also a $(1,0)$-vector field.
In fact, let $(z^1,\cdots, z^{n+l})$ be holomorphic coordinates. Then $X=\sum_{i=1}^{n+l}X^{i}\frac{\p}{\p z^{i}}$ and $Y=\sum_{i=1}^{n+l}Y^{i}\frac{\p}{\p z^{i}}$. Here $X^{i},Y^{i}$ may not be holomorphic functions.  Then
$$
\nabla_{\bar X} {Y}=\sum_{i,j} \bar X^{i}\left(Y^{j} \nabla_{\frac{\p}{\p \bar z^{i}}} \tfrac{\p}{\p z^{j}} +  \tfrac{\p Y^{j} }{\p \bar z^{i}} \tfrac{\p}{\p z^{j}}\right)=\sum_{i,j} \bar X^{i} \tfrac{\p Y^{j} }{\p \bar z^{i}} \tfrac{\p}{\p z^{j}}.
$$
We want  find a smooth  holomorphic $(n,0)$-field $L\left( \wedge_{j=1}^{n} S_{j}\right)\bigwedge \left( \wedge_{\alpha\in R_{M^{+}}} S_{\alpha}\right)$ on $\tau^{-1}(\Delta)$, where $L$ is a smooth function on $\tau^{-1}(\Delta).$
By the Koszul formula,  for any $j\leq n$, $\alpha\in R_{M^{+}}$, we have
\begin{align}
& g(\nabla_{\bar S_{j}}S_{\alpha},\bar S_{\alpha}) \label{eqn_Ch_4.1}\\
&=\frac{1}{2}\left[\bar S_{j}g(S_{\alpha},\bar S_{\alpha})-g(\bar S_{j},[S_{\alpha},\bar S_{\alpha}])+g(S_{\alpha},[\bar S_{\alpha},\bar S_{j}])+g(\bar S_{\alpha},[\bar S_{j},S_{\alpha}])\right] \nonumber\\
&=\frac{1}{2}\left[\frac{1}{4}\frac{\p D_{\alpha}}{\p x^{j}}-g(\bar S_{j},-\sqrt{-1}\sum_{k=1}^{n} M_{\alpha}^{k}H_{k})\right]+\frac{1}{4}D_{\alpha}\alpha(\tilde h_{j})=\frac{1}{4}D_{\alpha}\alpha(\tilde h_{j}).\nonumber
\end{align}
Similarly we have for any $\alpha,\beta\in R_{M^{+}},\;j,k,l\leq n,\;\;\;\;$
\begin{equation}\label{eqn_Ch_4.2}
g(\nabla_{\bar S_{j}}S_{l},\bar S_{k})=g(\nabla_{\bar S_{\alpha}}S_{l},\bar S_{k})=g(\nabla_{\bar S_{\alpha}}S_{\beta},\bar S_{\beta})=0.
\end{equation}
Denote $\mathbb{T}=\left( \wedge_{j=1}^{n} S_{j}\right)\bigwedge \left( \wedge_{\alpha\in R_{M^{+}}} S_{\alpha}\right)$.
By \eqref{eqn_Ch_4.1} and \eqref{eqn_Ch_4.2}, it is easy to see that
$$\nabla_{\bar S_{\beta}}\mathbb{T} = 0,\;\;\;\nabla_{\bar S_{j}}\mathbb{T} = \left(\sum_{\alpha\in R_{M^{+}}} \frac{  \alpha(\tilde h_{j})}{2} \right)\mathbb{T}.$$
Denote $\sigma_{i}=2\sum_{\alpha\in R_{M^{+}}}   { \alpha(\tilde h_{i}}).$
Choose $L=e^{-\frac{\sum_{i=1}^{n} \sigma_{i}x^{i}}{2}}$,  then
$$
\nabla_{\bar S_{\beta}}(L\mathbb T)=0,\;\;\;\nabla_{\bar S_{i}}(L\mathbb T)=0,\;\;\;\;\;\;i\leq n,\;\alpha\in R_{M^{+}}.
$$
  It follows that
$$
e^{\sum_{i=1}^{n}  \sigma_{i}x^{i}}\left( \wedge_{j=1}^{n} \omega^{j}\right)\bigwedge \left( \wedge^{\alpha\in R_{M^{+}}} \omega_{\alpha}\right)=|h|^{2}(dz^{i}\wedge d \bar z^{i} )^{n+l},
$$
for some holomorphic function $h.$ Then
we have
$$
{\Omega^{n+l}}=\det\left(\tfrac{\p^2 f}{\p x^{i}\p x^{j}}\right)\mathbb  D e^{-\sum_{i=1}^{n}  \sigma_{i}x^{i}} |h|^2{(\sum dz^{i}\wedge d\bar z^{i})^{n+l}}.
$$
Denote $\mathbb F_{\Delta}=\det(f_{ij})\mathbb D.$
Choose $\{S_{i},\bar S_{i},S_{\alpha},\bar S_{\alpha},i\leq n,\alpha\in R_{M^{+}}\}$ as a local frame field.  The Ricci tensor can be written as
 \begin{align}\label{eqn_Ric_4.1}
  \sum_{A, B} Ric(\frac{\p}{\p z^{A}},\frac{\p}{\p \bar z^{B}}) dz^{A}\wedge d\bar z^{B} =-\sum_{A,B}  \left[log \mathbb F_{\Delta}-\sum_{i=1}^{n}  \sigma_{i}x^{i}\right]_{,A\bar B}   \omega^{A} \wedge \bar \omega^{B}.
 \end{align}
 where $``,"$ denotes the covariant derivatives of $\mathcal G_{f}$ with respect to this frame field.

Since $\nabla_{\bar S_{A}}S_{B}$ is $(1,0)$-vector field and $\mathcal G( S_{\beta},\bar S_{l} )=0$, by Koszul formula, we have
\begin{align}
&\mathcal G(\nabla_{\bar S_{\alpha}}S_{\beta},\tfrac{\p}{\p x^{l}} )=\mathcal G(\nabla_{\bar S_{\alpha}}S_{\beta},\bar S_{l} ) \nonumber\\
&=\tfrac{1}{2}[-\bar S_{l}\mathcal G( \bar S_{\alpha}, S_{\beta})-\mathcal G(\bar S_{\alpha}, [S_{\beta},\bar S_{l}]) +\mathcal G(  S_{\beta}, [\bar S_{l},\bar S_{\alpha}])+\mathcal G(\bar S_{l}, [\bar S_{\alpha}, S_{\beta}])  ] \nonumber\\
&=\tfrac{1}{2}\delta_{\alpha\beta}\left[-\frac{1}{4}\frac{\p D_{\alpha}}{\p x^{l}}-\sum_{j}\frac{f_{jl}}{2}M_{\alpha}^{j}\right]=-\frac{1}{4}\frac{\p D_{\alpha}}{\p x^{l}}\delta_{\alpha\beta}.\label{eqn_Ch_4.3}
\end{align}
where we used $\frac{\p D_{\alpha}}{\p x^{l}}=\sum \frac{\p D_{\alpha}}{\p \xi_{k}}\frac{\p \xi_{k}}{\p x^{l}}=2\sum_{k} M_{\alpha}^kf_{kl}.$  Similarly we have
\begin{equation}\label{eqn_Ch_4.4}
\mathcal G(\nabla_{\bar S_{j}}S_{k},\tfrac{\p}{\p x^{l}} )=\mathcal G(\nabla_{\bar S_{j}}S_{\alpha},\tfrac{\p}{\p x^{l}} )=\mathcal G(\nabla_{\bar S_{\alpha}}S_{j},\tfrac{\p}{\p x^{l}} )=0.
\end{equation}
For any function  $F$ depending only on $(x^{1},\cdots, x^{n})$, we have
$$
F_{,A\bar B}=S_{A}\bar S_{B}F - \mathcal G(\nabla_{S_{A}}\bar S_{B},\tfrac{\p }{\p x^{l}})f^{kl}\tfrac{\p F}{\p x^{k}}.
$$
Then using  \eqref{eqn_Ch_4.3} and \eqref{eqn_Ch_4.4}  we have
\begin{align}\label{eqn_Co_4.5}
   F_{,j\bar k}=\tfrac{1}{4}\tfrac{\p^2 F}{\p x^{j}\p x^{k} },\;\; F_{,\alpha \bar j} =0,\;\;
 F_{,j\bar \alpha}=0\;\;\; F_{,\alpha \bar \beta}= \delta_{\alpha\beta}\tfrac{1}{4}\sum f^{kl}\tfrac{\p D_{\alpha}}{\p x^{k}}\tfrac{\p  F}{\p x^{l} }.
\end{align}
By \eqref{eqn_Ric_4.1} and \eqref{eqn_Co_4.5}, the Ricci curvatures are given by
\begin{align}\label{eqn_Ric_4.7}
&Ric (S_{j},\bar S_{k})=-\tfrac{1}{4}\tfrac{\p^2 log\mathbb F_{\Delta}}{\p x^{j}\p x^{k} },\;\;\;Ric (S_{\alpha},\bar S_{k})=
Ric  (S_{j},\bar S_{\alpha})=0,\;\;\; \\\label{eqn_Ric_4.8} & Ric (S_{\alpha},\bar S_{\beta})= \delta_{\alpha\beta}\left[ -\tfrac{1}{4}\sum f^{kl}\tfrac{\p D_{\alpha}}{\p x^{k}}\tfrac{\p log \mathbb F_{\Delta}}{\p x^{l} }+ \tfrac{1}{4} \sum \tfrac{\p D_{\alpha}}{\p \xi_{k}}\sigma_{k}\right],
\end{align}
The scalar curvature can be written as
\begin{equation}\label{eqn S_4.9}
 \mathbb{S}=-\sum_{i,j} f^{ij}\tfrac{\p^2 log\mathbb F_{\Delta}}{\p x^{i}\p x^{j} }- \sum_{k,l} f^{kl}\tfrac{\p \log\mathbb D }{\p x^{k}}\tfrac{\p log\mathbb F_{\Delta}}{\p x^{l} }+  h_G.
\end{equation}
In terms of $\xi$ and $u(\xi)$, $\mathbb{S}$ can be written as
\begin{equation}\label{eqn 4.6}
\mathbb{S}= -\frac{1}{\mathbb {D}}\sum_{i,j=1}^n\frac{\partial^2 \mathbb {D}u^{ij}}{\partial \xi_i\partial \xi_j} + h_G.
\end{equation}
Here
$$\mathbb {D}= \prod_{\alpha\in R_{M^{+}}} D_{\alpha},\;\;\;h_G= \sum \sigma_{i}\tfrac{\partial \log \mathbb{D}}{\partial \xi_{i}}
=\sum_{i,j}f^{ij} \sigma_{i}\tfrac{\partial \log \mathbb{D} }{\partial x^{j}},$$
where $\sigma $ is the sum of the positive roots of $R^+_M,$ and
$$\sigma_i=-2\sum_{\alpha\in R^{+}_M} \alpha(\sqrt{-1}  {\tilde H}_{i}).$$
Put $A=\mathbb{S}-h_G$. We will consider the PDE
\begin{equation}\label{eqn 4.7}
-\frac{1}{\mathbb {D}}\sum_{i,j=1}^n\frac{\partial^2 \mathbb {D}u^{ij}}{\partial \xi_i\partial \xi_j}=A.
\end{equation}
The equation \eqref{eqn 4.7} was introduced by Donaldson \cite{D5}
in the study of the scalar curvature of toric fibration, see also \cite{R} and \cite{N-1}. We call \eqref{eqn 4.7} a generalized Abreu Equation.

\v

 \section{Uniform stability}\label{Sec-Determinants}

We introduce several classes of functions. Set
\begin{align*}
\mc C&=\{u\in C(\bar\Delta):\, \text{$u$ is
convex on $\bar\Delta$ and smooth on $\Delta$}\},\\
\mathbf{S}&=\{u\in C(\bar\Delta):\, \text{$u$ is convex  on $\bar\Delta$
and $u-v$ is smooth on $\bar\Delta$}\},\end{align*}
where $v$ is given in \eqref{eqn2.1}.
For a fixed
point $p_o\in \Delta$, we consider
\begin{align*}
{\mc C}_{p_o}&=\{u\in \mc C:\, u\geq u(p_o)=0\},\\
\mathbf{S}_{p_o}&=\{ u\in \mathbf{S} :\, u\geq u(p_o)=0\}.\end{align*}
We say functions in ${\mc C}_{p_o}$ and ${\mathbf{S}}_{p_o}$ are {\it normalized} at $p_o$.
\v
We consider the generalized Abreu equation \eqref{eqn 4.7},
where, $\mathbb {D}>0$ and $A$ are given smooth functions on $\bar{\Delta}$.
 Following \cite{N-1} we consider the  functional
\begin{equation} 
\mc F_A(u)=-\int_\Delta \log\det(u_{ij})\mathbb {D}d \mu+\mc L_A(u),
\end{equation}
where
\begin{equation} 
\mc L_A(u)=\int_{\partial\Delta}u \mathbb {D}d\sigma-\int_\Delta Au \mathbb {D} d\mu.
\end{equation}
$\mc F_A$ is called the Mabuchi functional
and $\mc L_A$ is closely related to the Futaki invariants. The Euler-Lagrangian equation for $\mc F_A$ is \eqref{eqn 4.7}.
It is known that,
if $u\in \mathbf{S}$ satisfies the equation \eqref{eqn 4.7}, then $u$ is an absolute minimizer for
$\mc F_A$ on $\mathbf{S}$.

\begin{defn}\label{defn_1.5}
Let $\mathbb {D}>0$ and $A$ be smooth functions on $\bar\Delta$.
Then, $({\Delta},\mathbb {D},A)$ is called {\em uniformly $K$-stable} if
the functional $\mc L_A$ vanishes on affine-linear functions and
there exists a constant $\lambda>0$
such that, for any $u\in  {\mc C}_{p_o}$,
\begin{equation}\label{eqn 1.6}
\mc L_A(u)\geq \lambda\int_{\partial \Delta} u \mathbb {D}d \sigma.
\end{equation}
We also say that $\Delta$ is
$(\mathbb {D}, A,\lambda)$-stable.
\end{defn}
Using the same method in \cite{CLS4} we immediately get
\begin{theorem}\label{theorem_1.7}
If the equation \eqref{eqn 4.7} has a solution in $\mathbf{S}$, then $(\Delta, \mathbb {D}, A)$ is uniform K-stable.
\end{theorem}

\v \n
Donaldson \cite{D4} derived a $L^{\infty}$ estimate for the Abreu's equation in $\Delta \subset \mathbb {R}^2$. His method can be applied directly to the
generalized Abreu's Equation $\Delta \subset \mathbb {R}^2$ ( see also \cite{N-1}). We have
\begin{theorem}\label{theorem_3.2}
Let $\Delta\subset \mathbb {R}^2$ be a Delzant polytope, $\mathbb {D}>0$ and $A$ be two smooth functions defined on $\bar\Delta$. Let $u\in C^{\infty}(\Delta)$ satisfying \eqref{eqn 4.7}. Suppose that $\Delta$ is $(\mathbb {D},A,\lambda)$-stable. Then there is a constant
$\mff C_o>0$, depending on $\lambda$, $\Delta$, $\mathbb {D}$ and $\|A\|_{C^0}$, such that $|\max\limits_{\bar \Delta} u-\min\limits_{\bar \Delta} u|\leq \mff C_o$.
\end{theorem}
\v

\section{Differential inequalities}\label{Sec-differential inequalities}
\v
We first introduce some notations. Let $p$ be a vertex of $\Delta$, the edges meeting $p$ are the form $p+tE_i$, $t\geq 0$, $E^i\in \mathbb{Z}^n$, $i=1,...,n$. Consider the base $\{p; E^1,...,E^n\}$, let $\xi^p=(\xi_1^p,...,\xi_n^p)$ be the coordinates with respect to the base $\{p; E^i\}$. Suppose that
$$p=\sum_{i=1}^n c_i\nu^i,\;\; \nu^i=\sum_{i=1}^n a^i_j E^j,\;\;(a_{i}^j)\in SL(n,\mathbb Z),$$
we have the coordinate transformation
$$\xi^{p}_{i}=\sum_{j=1}^n a_{i}^j(\xi_{j}-c_{i}),\;\;\xi^{p}(p)=0,\;\;\;\Delta\subset \{\xi^{p}|\xi^{p}_{i}> 0 \}.$$
Set
$$x^{i}_{p}=\frac{\p u}{\p \xi_{i}^{p}},\;\;f_p(x^1_p,...,x^n_p)=\sum x^i_p \xi^p_i - u.$$

Then  $f_{p}$  can be naturally extend to a smooth function on a neighborhood $U_p$ of $\tau_{M}^{-1}(p).$ Let $$w^{i}=x^{i}_{p}+\sqrt{-1}y^{i}_{p},\;\;\;z^{i}_{p}=e^{\frac{w^{i}_{p}}{2}},\;\;\;i=1,\cdots, n .$$ Then $(z^{1}_{p},\cdots, z^{n}_{p})$ is a local holomorphic coordinates of $U_{p}.$
We have (see \cite{CLS2})
\begin{equation} \label{eqn_f_5.1}
f_{p}=f-\sum_{i=1}^n c_{i}x^{i}.
\end{equation}
Denote
\begin{align} \label{eqn_f_5.2}
&\mathbb F_{p}=4^{2n}\det\left(\tfrac{\p ^2 f_p}{\p z^i_{p} \p \bar  z^j_{p}}\right)\mathbb D =\det\left(\tfrac{\p ^2 f_p}{\p x^i_{p} \p x^j_{p}}\right)e^{-x^1_{p}-x^{2}_{p}-...-x^n_p}\mathbb D .
\\ \label{eqn_f_5.3}
&\mathbf{\Psi}_p:=\|\nabla log \mathbb F_{p}\|^2_f,\;\;\;
P=\exp\left(\kappa
\mathbb  F_p^{a}\right)\sqrt{\mathbb  F_p}{\mathbf{\Psi}_p},
\end{align}
where $a$ and $\kappa$ are positive constants to be determined later.
\v
More general, for each $(n-k)$-dimensional face of $\Delta$, one can associate it a
complex coordinate chart  of $M$:
$$\CHART_{k}\cong \cplane^{k}\times \left(\cplane^\ast\right)^{n-k}.$$
Let $(w^{1},\cdots, w^{n})$, $w^{i}=x^i + \sqrt{-1}y^i$, be the  log-affine coordinate of $M$. Set
$$z^{i}=e^{\frac{w^{i}}{2}},\;\;\;i\leq k.$$ Then
$(z^1,...,z^{k},w^{k+1},...,w^n)$ is local holomorphic coordinates of $\CHART_{k}.$
Denote $$Z_{n-k}=\{p\in M|z^{i}=0,i\leq k\},\;\;\;\;\;E=\{\xi|\xi_{i}=0,\;i\leq k\},$$
where $\xi_{i}=\frac{\p f}{\p x^{i}}$.
 Let $f_{E}=f-\sum_{i=1}^k c_{i}x^{i}-d,$  where $c_{i}$, $d$ are constants  such that $\tau_{f_{E}}(Z_{n-k})=E$.
Then   $f_{E}$  can be naturally extend to a smooth function in $\CHART_{k}$.
Denote
\begin{align} \label{eqn_f_5.2}
&\mathbb F_{E} =\det\left(\tfrac{\p ^2 f_E}{\p x^i  \p x^j }\right)e^{-\sum_{i=1}^{k} x^{i}}\mathbb D .
\\ \label{eqn_f_5.3}
&\mathbf{\Psi}_E:=\|\nabla log \mathbb F_{E}\|^2_f,\;\;\;
P_{E}=\exp\left(\kappa
\mathbb  F_E^{a}\right)\sqrt{\mathbb  F_E}{\mathbf{\Psi}_E}.
\end{align}
Set $V_{E}=\log \mathbb F_{E}$, $\mathbb A_{E}=-\square \log \mathbb F_{E}.$
\v
Let $\Omega_g$ be the Guillemin metric on $G\times_{K}M$ with local potential function $g$. Denote by $
\dot{R}_{i{\bar{j}}k\bar{l}}$ and $ \dot{R}_{i\bar{j}}$ the
curvature tensor and the Ricci curvature of $\Omega_g$,
respectively. Put
$$ \mc {\dot
{R}}:= \sqrt{\sum g^{m\bar{n}}g^{k\bar{l}}g^{i\bar{j}}g^{s\bar{t}}
 \dot{R}_{{m}\bar{l}i\bar{t}} \dot{R}_{k{\bar{n}}s\bar{j}}}.$$
Let $\Omega_f$ be a K\"ahler metric on $G\times_{K}M$ in the same K\"ahler class as $\Omega_g$  with local potential function $f$. Then there is a global defined function $\phi \in C^{\infty}(G\times_{K}M)$ such that
$\phi=f-g.$ Denote $\square =\sum f^{A\bar B} \frac{\p^2}{\p z^{A}\p \bar z^{B}}$ the Laplacian operator on $G\times_KM$ with respect to the metric $\Omega_{f}.$ Obviously $
n-\square \phi=\sum f^{i\bar j}g_{i\bar j}:=T.$  Put
\begin{equation*}
  Q=e^{-\mff N_1(\phi-\inf\phi+1)}\sqrt{\mathbb F_p}T,\;\;\;\;
\end{equation*} where $\mff N_1>0$ is a constant.
\v
In this section we establish some differential inequalities in $G(U_p)$ for $\log \mathbb{F}$, $P$ and $Q$.

\subsection{\bf Subharmonic function $\log \mathbb F_{p}+Nf_{p}$}
\v

 For any function $F$ depending only on $(x^{1},\cdots, x^{n})$
 we have
\begin{equation}\label{eqn_sq_5.2}
\square  F=\sum f^{ij}F_{ij} +\sum f^{ij}\frac{\log \mathbb D}{\p x^{i}} F_{j}=\sum f^{ij} F_{ij} +\sum \frac{\log \mathbb D}{\p \xi_{j}} F_{j}.
\end{equation}
 where $F_{i}=\frac{\p F }{\p x^{i}},F_{ij}=\frac{\p^2 F }{\p x^{i} \p x^{j}}$.
\begin{lemma}\label{lemma_5.1}

Choose $\{o,\nu^{i},\;i=1,...,n\}$ as a base of $\mathbb{R}^n$, let $\bar{\Delta}\subset \{(\xi_1,...,\xi_n)| \xi_1>0,\;\xi_2>0,\;...,\xi_n>0\}$ be a Delzant polytope satisfying
\begin{equation}\label{equ_R_3.6}
 \sum_{\alpha\in R_{M^{+}}}\frac{\sum_{j=1}^{n} M_{\alpha}^{j}diam(\Delta)}{D_{\alpha}}<\frac{n}{4}.
\end{equation}
Then there is a constant $N>0$ depending only on $n$, $\mathbb{D}$, $\Delta$ and the position of $\Delta$ in $\mathbb{R}^n$ such that for any vertex $p$ of $\Delta$
$$
\square (\log \mathbb F_{p}+Nf_{p})>0.
$$
\end{lemma}
\v\n
{\bf Proof.}  By \eqref{eqn_f_5.1} we have
 \begin{align}
\square  f_{p}&= \square  (f-\sum c_{i}x^{i}) =n +\sum \frac{\p \log \mathbb D}{\p \xi_{i}} \left(\frac{\p f}{\p x^{i}}-c_{i}\right) \nonumber\\
&  \geq n-\sum_{\alpha\in R_{M^{+}}}\frac{\sum_{j}M_{\alpha}^{j} diam(\Delta)}{\sum_{j}M_{\alpha}^{j}\xi_{j}}\geq \frac{n}{2}.
\end{align}
As $x^{i}=\frac{\p u}{\p \xi_{i}}=\sum_{j}\frac{\p u}{\p \xi^{p}_{j} }\frac {\p \xi^{p}_{j} } {\p \xi_{i} }=\sum a_{i}^jx_{p}^{j},$ we have $\det\left(\tfrac{\p ^2 f}{\p x^i_{p} \p x^j_{p}}\right)=\det\left(\tfrac{\p ^2 f}{\p x^i  \p x^j }\right).$
Denote by $(b_{i}^j)$ the inverse of $(a_{i}^j).$
Then by \eqref{eqn 4.6} and \eqref{eqn_sq_5.2}
\begin{align}
\square \log \mathbb  F_{p} =-\mathbb{S}+ \sum_{k=1}^{m}\tfrac{\p \log \mathbb  D}{\p \xi^{p}_{k}}(\sigma_{k}-\sum_{j} b^{j}_{k}):=-\mathbb  A_{p}.
\end{align}
Denote $\mathcal{P}$ be the set of all vertices of $\bar \Delta.$
Let $N$ be the constant such that
$$
\max_{p\in \mathcal{P}}\max_{\xi\in \bar\Delta} \mathbb A_{p}<\frac{n}{2}N.
$$
Then for any $p\in \mathcal{P},$ we have
$$
\square (\log \mathbb F_{p}+Nf_{p})>0.
$$
$\blacksquare$
\v

\subsection{\bf Differential inequalities for $P$}

\v
We are
going to calculate  $\df P$ and derive a differential inequality.
This inequality is similar to Lemma 5.1 in \cite{CLS1}. In this subsection we denote  $\mathbb F_{E},\mathbb A_{E},\mathbf \Psi_{E},\cdots$ by $\mathbb F,\mathbb A,\mathbf \Psi,\cdots,etc.$ By \eqref{eqn_Co_4.5}, \eqref{eqn_Ric_4.7} and  \eqref{eqn_Ric_4.8} we have
\begin{align}
\label{eqn_F_4.18} & (\log \mathbb  F )_{,j \bar k}=-Ric(S_{j},\bar S_{k}), \;\;\;(\log \mathbb  F )_{,\alpha \bar k} =0, \;\;\; (\log \mathbb  F )_{,k\bar \alpha }=0, \;\;\;
\\\label{eqn_F_4.19}
&
(\log \mathbb  F)_{,\alpha \bar \beta}= -Ric(S_{\alpha},\bar S_{\beta})+\frac{\delta_{\alpha\beta}}{4}h_{\alpha},
\end{align}
where $h_{\alpha}=\sum_{k} \frac{\p D_{\alpha}}{\p \xi_{k}}\left(\sigma_{k}-\sum_{j} c^{j}_{k}\right)  $ and $c^{j}_{k}$ are constants depending only on $E.$
 Denote by $f_{,A\bar B}$ the components of the metric $\mathcal G$ with respect to the frame $\{S_A, \bar{S}_B\}$, and by $(f^{A\bar{B}})$ the inverse of $(f_{,A\bar B})$. Put $W=\det(f_{,i\bar j})$,  $ V=\log \mathbb  F $,
$$  \|V_{,A\bar B}\|_f^2=\sum f^{A\bar B}f^{C\bar D}V_{A\bar D}V_{B\bar C},\;\;
\|V_{,AB}\|_f^2=\sum f^{A\bar B}f^{C\bar D}V_{,AC}V_{,\bar B\bar D}. $$

\begin{lemma}\label{lemma_5.1.1}
\begin{equation*}
\frac{\df  P}{P}  \geq   \frac{\|V_{,A\bar{B}}\|_f^2}{2\mathbf{\Psi}}
 +a^2\kappa(1 -2 \kappa \mathbb  F^{a}  )\mathbb  F ^{a} \mathbf{\Psi}
-   \frac{2|\langle \nabla \mathbb A, \nabla V\rangle|}{\mathbf{\Psi}}- \left({a} \kappa \mathbb  F^{a}
+\tfrac{1}{2}\right)\mathbb  A,
\end{equation*}
where $\langle, \rangle$ denotes the inner product with respect to the metric $\Omega_f.$
\end{lemma}
\v\n
{\bf Proof.} By definition,
 \begin{equation*} \df {\mathbf{\Psi}} =\sum
f^{A\bar{B}}f^{C\bar{D}} \left(V_{,A} V_{,\bar{B}C\bar{D}} +
V_{,AC\bar{D}}V_{,\bar{B}} + V_{,AC}V_{,\bar{B}\bar{D}} +
V_{,A\bar{D}}V_{,\bar{B}C}\right).
 \end{equation*}
Since $V$ depends only on $(x^{1},\cdots, x^{n})$, we have
\begin{equation}\label{eqn_V_alp_1}
V_{,\alpha}=V_{,\bar \alpha}=0.
\end{equation}
  By the Ricci
identities, \eqref{eqn_V_alp_1} and $f^{j\bar \alpha}=0,$ we have
\begin{equation*}
V_{,\bar{B}C\bar{D}} = V_{,C\bar{D}\bar{B}},\;\;\; V_{,AC\bar{D}}
= V_{,C\bar{D}A}+\sum f^{j\bar{h}} V_{,j} R_{C\bar{h}A\bar{D}}.
\end{equation*}
It follows that
\begin{equation}\label{eqn_5.5}
\df {\mathbf{\Psi}}  =  \|V_{,AB}\|^2+\|V_{A\bar B}\|^2+
\sum f^{i\bar{j}}f^{k\bar{l}}\left(- V_{i\bar{l}}V_{,k}V_{,\bar{j}}
 \right) -
  2Re (\sum f^{i\bar{j}}V_{,i}\mathbb A_{,\bar{j}}),
   \end{equation}
  where we use the facts
$ R_{i\bar{j}}= - V_{k\bar{l}}$ and  $\df V=-\mathbb  A . $
Denote $\Pi=a \kappa \mathbb  F^{a} + \tfrac{1}{2}.
$ Then
\begin{eqnarray}\label{eqn_5.6}
 P_{,i}&=&  P\left(\frac{{\mathbf{\Psi}}_{,i}}{\mathbf{\Psi}} +
\Pi V_{,i}\right) =:   P\Lambda_i, \nonumber\\
\df   P&=&  P \left[\sum f^{i\bar j}\Lambda_i\Lambda_{\bar j}
+\frac{\df\mathbf{\Psi}}{\mathbf{\Psi}} -\frac{\|\nabla \mathbf{\Psi}\|_f^2}{\mathbf{\Psi}^2}
+\Pi\df V + a^2\kappa  \mathbb  F^{a} {\mathbf{\Psi}} \right].\;\;\;\;\;\;
\end{eqnarray}
For any point $q\in G(\mff U_{k}),$ we choose an affine transformation of the frame fields $\{S_{A},\bar S_{B}\}$      such
that, at $q$,
\begin{equation*}
f_{,i\bar{j}}= c\delta_{ij},\;\;1\leq i,j\leq n,\;\;\;\;\;\;\;\; V_1 = V_{\bar{1}},\;\;\;V_i =
V_{\bar{i}}=0,\;\;\forall \;1\leq i \leq n,\end{equation*} where
$c=[W(q)]^{\frac{1}{n}}.$ Then by the same arguments of Lemma 5.1 in \cite{CLS1}
 we can prove the lemma. $\blacksquare$ \v

\subsection{\bf Differential inequality for $Q$}\label{sect_5.3}

\v

 The following
differential inequality of $ n-\square \phi$ has been proved in \cite{CLS1} (see Section \S 5, Inequality-II)
\begin{lemma}\label{lemma_5.3}
\begin{equation}\label{eqn_5.9}
 \df   \log \left(n-\square \phi\right)  \geq  -\| Ric\|_f - \mathcal {\dot
R}(n-\square \phi).
\end{equation}
\end{lemma}
\v
In the following we restrict ourself to $n=2$.

\v
\begin{lemma}\label{lemma_5.3.1}
 Suppose that
  \begin{equation}\label{eqn_5.3.1}
  \|\mathbb A_{p}\|_{C^{1}(\tau_{f}(U_{p}))}\leq \mff N_2,\;\;\;\;\;\max_{\bar U_{p}} \mathbb F\leq\mff
N_2,\;\;\;\;\;\; \max_{\bar U_{p}}|\phi|+|z|\leq \mff N_2\end{equation}  for some constant
$\mff N_2>0.$ Then we may choose
\begin{equation}\label{eqn_5.3.2}
\mff N_1=100,a=\frac{1}{3},\kappa=[4\mff N_2^\frac{1}{3}]\inv
\end{equation}
 such that
\begin{equation}\label{eqn_5.3.3}
\df (P+Q+ \mff C_1f_p ) \geq \mff C_{2}(P+Q)^2>0
\end{equation} for some positive constants $\mff C_1$ and $\mff C_2$ that
depend only on $\mff N_2$, the structure constants of $\mathfrak{g}$, $\mathbb{D}$, $\Delta$ and the position of $\Delta$ in $\mathbb{R}^2$.
\end{lemma}
{\bf Proof.} Applying  Lemma \ref{lemma_5.1.1} and the choice of $a=\frac{1}{3}$ and $\kappa$, in particular,
$\kappa \mathbb F^a\leq 1/4$, we have
\def \al{\frac{1}{3}}
\begin{equation}\label{eqn_5.3.4}
\frac{\mathbf{\Psi}\df  P}{P} \geq \left(\frac{1}{2}\|V_{,A\bar{B}}\|_f^2
 +\frac{1}{18} \kappa \mathbb F^{\al}\mathbf{\Psi}^2\right)
- \left(2 |\langle \nabla \mathbb A, \nabla V\rangle| +\mathbf{ \Psi}
|\mathbb A| \right).
\end{equation}
Treatment for $\langle \nabla \mathbb A, \nabla V\rangle$:
using log-affine coordinates we have
\begin{equation*}
|\langle \nabla \mathbb A, \nabla V\rangle| =
\left|\sum f^{ij}\frac{\partial \mathbb A}{\partial x^{i}}\frac{\partial
V}{\partial x^{j}} \right| =\left|\sum  \frac{\partial K}{\partial
\xi_k} \frac{\partial V}{\partial
x^{k}} \right| \leq \mff N_2\sum_{j}\left|\frac{\partial V}{\partial x^{j}}
\right|.
\end{equation*}
If we use the complex coordinates $z_i$, we have
$$
\left|\frac{\partial V}{\partial x^{j}} \right|=   \left|z_j \frac{\partial V}{\partial z_j}\right|.
$$
Since $|z|$ is bounded, in this coordinates we have
$$
C^{-1} \leq  g_{i\bar j} \leq C,\;\; \sum f^{i\bar i}\leq C \sum f^{i\bar j}g_{i\bar j}\leq C T.
$$
Then we conclude that
\begin{equation*}
|\langle \nabla \mathbb A, \nabla V\rangle|  \leq C\sqrt{2\mathbb FT\mathbf{\Psi}}.
\end{equation*}

We explain the last step: suppose that $0<\nu_1\leq\nu_2$ are the
eigenvalues of $(f_{i\b j})$, then
$$
\mathbf{\Psi}=\sum f^{i\bar j}V_iV_{\bar j} \geq \nu_2\inv
(|V_1|^2+|V_2|^2)\geq(\mathbb FT)\inv(|V_1|^2+|V_2|^2).
$$
Note that $$(e\mathbb F)^{-\half}\leq \mathbf{\Psi}/P=(\exp(\kappa \mathbb F^\alpha)
\mathbb F^\half)\inv \leq \mathbb F^{-\half},$$ \eqref{eqn_5.3.4} is then
transformed to be
\begin{equation}\label{eqn_5.3.5}
{\df  P} \geq \mathbb F^\half\left(\frac{1}{2}\|V_{,A\bar{B}}\|_f^2
 +\frac{1}{18} \kappa \mathbb F^{\al}\mathbf{\Psi}^2\right)
- C' \mathbb F^\half\left(\sqrt{\mathbb FT\mathbf{\Psi}} + \mathbf{\Psi} |\mathbb A| \right).
\end{equation}
Applying the Young inequality and the Schwartz inequality to  terms
in \eqref{eqn_5.3.5},
  we  have that
\begin{equation}\label{eqn_5.3.7}
\df  P \geq \half \mathbb F^\half \|V_{,A\bar B}\|_f^2 + C_1\mathbb F^{-\frac{1}{6}}
P^2  -\epsilon QT -C_2(\epsilon).
\end{equation}
By \eqref{eqn_F_4.18}, \eqref{eqn_F_4.19} and  Cauchy inequalities, for any $\delta\in (0,1),$ we obtain that
  $$\|Ric\|^2_{\mathcal G_f}=\sum_{i,j} \|V_{i\bar j}\|^2 +\sum_{\alpha\in R_{M^{+}}}\|V_{,\alpha\bar \alpha}+h_{\alpha}\|^2\leq (1+\delta)\|V_{,A\bar B}\|^2 +C_{{\delta} ,  R_{M^{+}}} ,$$
  where $C_{{\delta} ,  R_{M^{+}}}>0$ is a constant depending only on $1/\delta$ and $\sum_{\alpha\in R_{M^{+}}}h_{\alpha}^2.$
 By a direct calculation and  the formula \eqref{eqn_5.9} we have
\begin{align*}\df Q &\geq Q\left( N_1 T+\half \square V +
 \df \log T
\right) \\
 &\geq   Q \left((N_1-\mathcal {\dot R})T-\half\mathbb A-(1+\delta)\left\|V_{,A\bar{B}}\right\|_f-C_{{\delta} ,  R_{M^{+}}}\right),
\end{align*}
 Choosing $\delta $ small,  using the explicit value $N_1$ and the bounds of $\mathbb F$ and $\mathbb A$,
 applying the Schwartz inequality properly,
we can get
\begin{equation}\label{eqn_5.3.6}
\df Q\geq   -\frac{1}{4} \mathbb F^\half\|V_{,A\bar B}\|_f^2 + \frac{\mff N_1}{3}
QT -C_3(\mff N_1,\mff N_2,\delta,R_{M^{+}}).\end{equation}
Combining \eqref{eqn_5.3.7} and \eqref{eqn_5.3.6}, and choosing $\epsilon=\frac{1}{100}$, we have
\begin{equation*}
 {\df (  Q+  P)} \geq  {C_1\mathbb F^{-\six} P^2+\frac{N_1}{4}
QT}  - {C_4} .
\end{equation*}
 Note that
$$
T= e^{-N_1\phi}\mathbb F^{-\half}  Q\geq
e^{-N_1\phi}\mff N_2^{-\al}\mathbb F^{-\six} Q \geq C_5\mathbb F^{-\six}  Q,
$$
 we get
\begin{equation*}
{\df (  Q+  P)}\geq C_6\mathbb F^{-\six}(  Q+  P)^2 -C_7
\end{equation*}
where  $C_6,C_7$ are constants depending only on $C_1,\mff N_1$ and $\mff
N_2$. Our lemma follows from  $
 \Box  f_{p} \geq n/2.
 $  and $|\mathbb F|\leq \mff N_2$.
$\;\;\;\; \blacksquare$

\v

\subsection{Interior estimate of $\mathbf\Psi$}\label{sect_6.1}

\v

We use Lemma \ref{lemma_5.1.1} to
derive the interior  estimate of   $\mathbf{\Psi}$ in a geodesic ball in $\CHART_{k}$. In this subsection we denote  $\mathbb F_{E},\mathbb A_{E},\mathbf \Psi_{E},\cdots$ by $\mathbb F,\mathbb A,\mathbf \Psi,\cdots,etc.$

\begin{lemma}\label{lemma_7.16a} Let  $B_{a}(o)\subset \CHART_{k}$ be a closed
geodesic ball of radius $a$ centered at o with respect to the metric $\mathcal G_{f}$. Set $\mathbb F_\diamond:=\max\limits_{B_{a}(o)} \mathbb F.$  Suppose that
\begin{equation} \label{eqn_7.31a}
\min\limits_{B_{a}(o)}|\mathbb A|\neq0,
\mbox{ and, }  \mathbb F^{\half} (\mc K+\mathbb K+\|\nabla \log |\mathbb A|\|^2_f+ \mathbf{\Psi})\leq 4,
\end{equation} in $B_{a}(o)$. Then the following estimate  holds in
$\; B_{a/2}(o)\;$
\begin{eqnarray}\label{eqn_7.32a} {\mathbb F^{\half}}  \mathbf{\Psi}
 \leq \mff C_{6}\left[   \mathbb F_\diamond ^{\half}\max_{B_{a}(o)}|\mathbb A| +  \mathbb F_\diamond ^{\frac{1}{3}}\max_{B_{a}(o)}|\mathbb A|^{\frac{2}{3}}  +
a\inv \mathbb F_\diamond^{\frac{1}{4}}+ a^{-2}\mathbb F_\diamond ^{\half} \right],\;\;
   \;
 \end{eqnarray}    where $\mff C_6$ is a constant  depending
only on $n$.
 \end{lemma}
\v \n {\bf Proof.} Consider the function
$$F:= ( a^2-r^2)^2 {P} $$
defined in   $B_{a}(o),$ where $r$ denotes the geodesic distance from $o$ to $z$ with
respect to the metric $ \mathcal G_f.$  $F$ attains its
supremum at some interior point $q^\ast$.
Choose
$\kappa= \frac{1}{4 {\mathbb F_\diamond}^{\frac{1}{4}}},\;\;\;\;\alpha=\frac{1}{4},$ As in \cite{CLS1}, using Lemma \ref{lemma_5.1.1} and a direct calculation we have,  at $q^\ast$,
\begin{equation}\label{eqn_7.35a}
 \epsilon_{o}\left(\frac{\mathbb F}{\mathbb F_\diamond}\right)^{\tfrac{1}{4}} \mathbf{ \Psi}
-   \frac{2\|\nabla  {\mathbb A}\|_{f}}{\sqrt{\mathbf{\Psi}}}- \frac{9}{16} {|\mathbb A|}
 -\frac{24 {a}^2}{(a^2-r^2)^{2}}
 -\frac{4(1+r\df r)}{a^2-r^2}\leq 0,\end{equation}
  where $\epsilon_{o}=\frac{1}{128}.$
  Denote by $\Gamma$ the geodesic from $o$ to $q^\ast$.  To estimate $r\df r$, we consider two cases:
\v\n
{\bf Case 1.} Along $\Gamma$ the following estimate holds
$$\mathbb F\geq \frac{1}{100}\mathbb F(q^\ast).$$
{\bf Case 2.} There is a point $q\in \Gamma$ such that $\mathbb F(q)<\frac{1}{100}\mathbb F(q^\ast)$. Then there is a point $q_1\in \Gamma$ such that
$$\mathbb F(q_1)=\frac{1}{100}\mathbb F(q^\ast),\;\;\;\;\;\;\mathbb F(q)\geq \frac{1}{100}\mathbb F(q^\ast),\forall q\in [q_1,q^\ast].$$
Then by the same argument of \cite{CLS1} for both cases we have  estimates \eqref{eqn_7.32a} in $B_{a/2}(o)$. $\blacksquare$

\v\v\v

\bibliographystyle{elsarticle-num}
\bibliography{<your-bib-database>}

\end{document}